\input amstex
\documentstyle{amsppt}
\magnification1200 \NoBlackBoxes
\pagewidth{6.5 true in} \pageheight{9 true in}

\topmatter\title Primitive prime factors in second-order linear recurrence sequences \endtitle
\author Andrew Granville
\endauthor

\address{D\' epartement de Math\' ematiques et statistique,
Universit\' e de Montr\'eal, CP 6128 succ. Centre-Ville, Montr\' eal
QC  H3C 3J7, Canada } \endaddress
\email{andrew{\@}dms.umontreal.ca} \endemail


\thanks{Many thanks are due to Hendrik Lenstra for  informing me of his observation that one can determine the value of the Jacobi symbol $(x_m/x_n)$, where $x_m=2^m-1$,  in terms of the inverse of $m$ mod $n$. He also proved our results for odd $n$ and then we shared several further email exchanges that led to the results discussed in this paper.  I would also like to thank Ilan Vardi for the email discussion alluded to herein, and of course NSERC for partially funding this research.\newline
{\tt 2010 Mathematics Subject Classification}: 11B37. Secondary: 11A05, 11D41 \   \newline
{\tt Keywords}:  Recurrence sequences; primitive prime factors; characteristic prime factors \
 }\endthanks

\dedicatory{To Andrzej Schinzel on his 75th birthday, with thanks for the many inspiring papers}
\enddedicatory

\abstract For a class of Lucas sequences $\{ x_n\}$, we show that
if $n$ is a positive integer then $x_n$ has a primitive prime
factor which divides $x_n$ to an odd power, except perhaps when
$n=1, 2, 3$ or $6$. This has several desirable consequences.
\endabstract
\endtopmatter

\parindent=20pt

\document

\head 1. Introduction \endhead

\subhead 1a. Repunits and primitive prime factors \endsubhead
The numbers $11$, $111$ and $1111111111$ are known as {\sl
repunits}, that is all of their digits are 1 (in base $10$).
Repunits cannot be squares (since they are $\equiv 3 \pmod 4$), so
one might ask whether a product of distinct repunits can ever be a
square? We will prove that this cannot happen. A more interesting
example is the set of repunits in base 2, the integers of the form
$2^n-1$. In this case there is one easily found product of
distinct repunits that is a square, namely $(2^3-1) (2^6-1) =
21^2$ (which is  $111\cdot 111111 = 10101 \cdot 10101 $ in base
2); this turns out to be the only example.

For a given sequence of integers $\{ x_n\}_{n\geq 0}$, we define a
{\sl characteristic prime factor} of $x_n$ to be a prime $p$ which
divides $x_n$ but $p\nmid x_m$ for $1\leq m\leq n-1$. The
Bang-Zsigmondy theorem (1892) states that if $r>s\geq 1$ and
$(r,s)=1$ then the numbers
$$ x_n=\frac{r^n-s^n}{r-s} $$
have a characteristic prime factor for each $n>1$ except for the case $\frac{2^6-1}{2-1}$. A  {\sl primitive prime factor} of $x_n$ is a characteristic prime factor of $x_n$ that does not divide $r-s$.

For various Diophantine applications it would be of interest to
determine whether there is a characteristic prime factor $p$ of $x_n$
for which $p^2$ does not divide $x_n$.  As an example of such an
application, note that if $x_{n_1}\dots x_{n_k}$ is a square where
$1<n_1<n_2<\dots <n_k$ and $k\geq 1$ then a characteristic prime factor
$p$ of $x_{n_k}$ divides only $x_{n_k}$ in this product and hence
must divide $x_{n_k}$ to an even power. Thus if $p$ divides
$x_{n_k}$ to only the first power then $x_{n_1}\dots x_{n_k}$
cannot be a square. Unfortunately we are unable to prove anything
about characteristic prime factors dividing  $x_n$ only to the first
power, but we are able to  show that there is a characteristic prime
factor that divides $x_n$ to an odd power, which is just as good
for this particular application.

\proclaim{Theorem 1}   If $r$ and $s$ are pairwise coprime integers for which 2 divides $rs$ but not 4, then $(r^n-s^n)/(r-s)$ has a characteristic prime factor which divides it to an odd power, for each $n>1$ except perhaps for $n=2$ and $n=6$. The case $n=2$ is exceptional if and only if $r+s$ is a square. The case $n=6$ is exceptional if and only if $r^2-rs+s^2$ is 3 times a square.
\endproclaim

In particular $2^n-1$ has a characteristic prime factor which divides it to an odd power, for all $n>1$ except $n=6$. Also $(10^n-1)/9$ has a characteristic prime factor which divides it to an odd power for all $n>1$. One can take these all to be primitive prime factors.

\proclaim{Corollary 1} Let $x_n=(r^n-s^n)/(r-s)$ where $r$ and $s$ are pairwise coprime integers for which 2 divides $rs$ but not 4.  If $x_{n_1}x_{n_2} \dots x_{n_k}$ is a square where $1<n_1<n_2<\dots <n_k$ and $k\geq 1$, then either $x_2=r+s$ is a square, or  $x_3x_6=x_3^2(r^3+s^3)$  is a square.
\endproclaim

The infinitely many examples of this last case include $2^3+1=3^2$, leading to the solution $(2^3-1)(2^6-1) = 21^2$, and
$74^3-47^3=549^2$ leading to
 $\frac{74^3-(-47)^3}{121} \cdot
\frac{74^6-(-47)^6}{121} = 2309643^2$. Since $2^3+1=3^2$ is the
only non-trivial solution in integers to $r^3+1=t^2$, we have
proved that the only example of a product of repunits which equals
a square,  in any base $b$ with $b\equiv 2 \pmod 4$, is the one
base 2 example $(2^3-1)(2^6-1) = 21^2$ given already.

\bigskip

\subhead 1b. Certain Lucas sequences \endsubhead
The numbers $x_n=(r^n-s^n)/(r-s)$ satisfy $x_0=0, x_1=1$ and the
second order linear recurrence $x_{n+2}=(r+s)x_{n+1}-rsx_n$ for each $n\geq 0$. These are examples of a  Lucas sequence, where $\{ x_n\}_{n\geq 0}$ is a {\sl Lucas sequence} if  $x_0=0, x_1=1$ and
$$
x_{n+2}=bx_{n+1}+cx_n \ \ \text{ for all} \ \  n\geq 0 , \tag{1}
$$
for given non-zero, coprime integers $b,c$. The {\sl discriminant} of the Lucas sequence is
$$
\Delta:= b^2+4c.
$$
Carmichael showed in 1913 that if   $\Delta>0$  then $x_n$ has a characteristic prime factor for each $n\ne 1, 2$ or $6$ except for
$F_{12}=144$ where  $F_n$ is  the Fibonacci sequence ($b=c=1$), and for  $F'_{12}$ where $F'_n=(-1)^{n-1}F_n$ ($b=-1,\ c=1$). Schinzel [7], defined a {\sl primitive prime factor} of $x_n$ to be a characteristic prime factor of $x_n$ that does not divide the discriminant $\Delta$.

We have been able to show the analogy to Theorem 1 for a class of Lucas sequences:

\proclaim{Theorem 2} Let $b$ and $c$ be pairwise coprime integers with $c\equiv 2\pmod 4$  and $\Delta=b^2+4c>0$.
Let $\{ x_n\}_{n\geq 0}$ be the Lucas sequence satisfying (1). If $n\ne 1, 2$ or $6$ then
$x_n$ has a characteristic prime factor which (exactly) divides $x_n$ to an odd power.
\endproclaim

In fact $x_2$ does not have such a prime factor if and only if $x_2=b$ is a square; and $x_6$ does not have such a prime factor if and only if $x_6/(x_3x_2)=b^2+3c$ equals 3 times a square.

Theorem 1 is a special case of Theorem 2 since there we have
$c=-rs\equiv 2\pmod 4$, \ $(b,c)=(r+s,rs)=1$ and
$\Delta=(r-s)^2>0$.

\proclaim{Corollary 2} Let the Lucas sequence $\{ x_n\}_{n\geq 0}$
be as in Theorem 2.  If $x_{n_1}x_{n_2} \dots x_{n_k}$ is a square
where $1<n_1<n_2<\dots <n_k$ and $k\geq 1$ then the product is
either $x_2=b$ or $x_3x_6$.
\endproclaim

In fact $x_3x_6$ is a square if and only if $b$ and $b^2+3c$ are both 3 times a square; that is, there exist odd integers $B$ and $C$ with $(C,3B)=1$ and $4C^2>3B^4$, for which $b=3B^2$ and $c=C^2-3B^4$.

With a little more work we can improve Theorem 2 to account of the notion of primitive prime factors:

\proclaim{Theorem 3} Let $b$ and $c$ be pairwise coprime integers with $c\equiv 2\pmod 4$  and $\Delta=b^2+4c>0$.
Let $\{ x_n\}_{n\geq 0}$ be the Lucas sequence satisfying (1). If $n\ne 1, 2, 3$ or $6$ then
$x_n$ has a primitive prime factor which (exactly) divides $x_n$ to an odd power.
\endproclaim

The exceptions for $n=1,2$ and $6$ are as above in Theorem 2. In fact $x_3$ does not have such a prime factor if and only if $x_3=b^2+c$ equals 3 times  a square;

\subhead 1c.  Fermat's last theorem and  Catalan's conjecture; and a new observation \endsubhead
Before Wiles' work, one studied Fermat's last theorem by considering the equation $x^p+y^p=z^p$ for prime exponent $p$ where $(x,y,z)=1$, and split into two cases depending on whether $p$ divides $xyz$. In the ``first case'', in which $p\nmid xyz$, one can factor $z^p-y^p$ into two coprime factors $(z-y)$ and
$(z^p-y^p)/(z-y)$ which must both equal the $p$th power of an integer. Thus if the $p$th term of the Lucas sequence $x_p=(z^p-y^p)/(z-y)$ is never a $p$th power for odd primes $p$ then the first case of Fermat's last theorem follows, an approach that has not yet succeeded. However Terjanian [9] {\sl did} develop these ideas to prove that the first case of Fermat's last theorem is true for even exponents, showing that if $x^{2p}+y^{2p}=z^{2p}$ in coprime integers $x,y,z$ where $p$ is an odd prime then $2p$ divides either $x$ or $y$:

In any solution,  $x$ or $y$ is even, else $2$ divides $(x^{p})^2+(y^{p})^2=z^{2p}$ but not 4, which is impossible.
So we may assume that  $x$ is even, but not divisible by $p$, and $y$ and $z$ are odd so that we have a solution $r=z^2,\ s=y^2,\ t=x^p$ to $r^p-s^p=t^2$  with $r\equiv s\equiv 1 \pmod 4$ and $(t,2p)=2$.  Let $x_n=(r^n-s^n)/(r-s)$ for all $n\geq 1$,
so that $x_p(r-s)=t^2$ and $(x_p,r-s)=(p,r-s)|(p,t)=1$, which implies that $x_p$ is a square. Terjanian's key observation is that the Jacobi symbols
$$
\left(  \frac{x_m}{x_n} \right) = \left(  \frac{m}{n} \right) \ \ \text{for all odd, positive integers} \ m \ \text{and} \ n. \tag{2}
$$
Thus by selecting $m$ to be an odd quadratic non-residue mod $p$,
we have $(x_m/x_p)=-1$ and therefore $x_p$ cannot be a square.
This contradiction implies that $p$ must divide $t$, and hence Terjanian's result.

A similar method was used earlier by Chao Ko [2] in his proof
that $x^2-1=y^p$ with $p>3$ prime has no non-trivial solutions (a
first step on the route to proving Catalan's conjecture).
Rotkiewicz [4] showed, by these means, that if $x^p+y^p=z^2$ with
$(x,y)=1$ then either $2p$ divides $z$ or $(2p,z)=1$, which
implies both Terjanian's and Chao Ko's results. Rotkiewicz's key
lemma in [4], and then his Theorem 2 in [5], extend (2):\
Assume that $\Delta$ and $b$ are positive with $(b,c)=1$. If $b$
is even and $c\equiv -1 \pmod 4$ then (2) holds.  If 4 divides
$c$, or if $b$ is even and $c\equiv 1 \pmod 4$ then $(x_m/x_n)=1$
for all odd, positive integers $m,n$. In the most interesting
case, when $2$, but not $4$, divides $c$,  we have
$$
\left(  \frac{x_m}{x_n} \right) = (-1)^{\Lambda(m/n)} \ \ \text{for all odd, coprime, positive integers} \ m \ \text{and} \ n>1, \tag{3}
$$
where $\Lambda(m/n)$ is the length of the continued fraction for
$m/n$; more precisely, we have a unique representation $m/n=[a_0,a_1,\dots,a_{\Lambda(m/n)-1}]$ where
each $a_i$ is an integer, with $a_0\geq 0,\ a_i\geq 1$ for each $i\geq 1$, and $a_{\Lambda(m/n)-1}\geq 2$.

Note that we have not given an explicit evaluation of $(x_m/x_n)$
when $b$ and $c$ are both odd, the most interesting case being
$b=c=1$ which yields the Fibonacci numbers. Rotkiewicz [6] does give
a complicated formula for determining $(F_m/F_n)$ in terms of a special continued fraction type expansion for $m/n$; it remains to find a simple way to evaluate this formula.

To  apply (3) we show that one can replace $\Lambda(m/n) \pmod 2$ by the much simpler $[2u/n] \pmod 2$, where $u$ is any integer $\equiv 1/m \pmod n$ (and that this formula holds for all coprime positive integers $m, n$). Our proof of this, and the more general (4), is direct (see Theorem 4 and Corollary 6 below), though Vardi explained, in email correspondence, how to use the theory of continued fractions to show that  $\Lambda(m/n) \equiv [2u/n] \pmod 2$ (see the end of section 5).

It is much more difficult to prove that Lucas sequences with negative discriminant have primitive prime factors. Nonetheless, in 1974 Schinzel [8] succeeded in showing that $x_n$ has a primitive prime factor once $n>n_0$, for some sufficiently large $n_0$,
if  $\Delta\ne 0$, other than in the periodic case $b=\pm 1,\ c=-1$. Determining the smallest possible value of $n_0$ has required great efforts culminating in the beautiful work of Bilu, Hanrot and Voutier [1] who proved that $n_0=30$ is best possible. One can easily deduce from Siegel's theorem that if $\phi(n)>2$ then there are only finitely many Lucas sequences for which $x_n$ does not have a primitive prime factor, and these exceptional cases are all explicitly given in [1]. They show that such examples occur only for $n=5, 7, 8, 10, 12, 13, 18, 30$: if  $b=1, c=-2$ then $x_5, \ x_{8},\ x_{12},\ x_{13},\ x_{18},\ x_{30}$ have no primitive prime factors; if  $b=1, c=-5$ then $x_7=1$; if $b=2, c=-3$ then $x_{10}$ has no primitive prime factors; there are a handful of other examples besides, all with $n\leq 12$.

\subhead 1d. Sketches of some proofs \endsubhead
In this subsection we sketch the proof of a special case of Theorem 2 (the details will be proved in the next four sections). The reason we focus now on a special case is that this is already sufficiently complicated, and extending the proof to all cases involves some additional (and not particularly interesting) technicalities, which will be given in section 6.

\proclaim{Theorem 2'} Let $b$ and $c$ be integers for which $b\equiv 3 \pmod 4, \ c\equiv 2\pmod 4$, the Jacobi symbol $(c/b)=1$ and $\Delta=b^2+4c>0$. If $\{ x_n\}_{n\geq 0}$ is the Lucas sequence satisfying (1) then
$x_n$ has a characteristic prime factor which (exactly) divides $x_n$ to an odd power for all $n>1$ except perhaps when $n=6$.
This last case occurs if and only if $x_6/(3x_2x_3)$ is a square.
\endproclaim

\demo{Sketch of the proof of Theorem 2'}
Let $x_n=y_nz_n$ where $y_n$ is divisible only by characteristic prime factors of $x_n$, and $z_n$ is divisible only by non-characteristic prime factors of $x_n$.  If every characteristic prime factor divides $x_n$ to an even power then $y_n$ is a square: it is our goal to show that this is impossible.

A complex number $\xi$ is a {\sl primitive $n$th root of unity}  if $\xi^n=1$ but $\xi^m\ne 1$ for all $1\leq m< n$. Let $\phi_n(t)\in \Bbb Z[t]$ be the {\sl $n$th cyclotomic polynomial}, that is the monic polynomial whose roots are the primitive $n$th roots of unity. Evidently $x^n-1=\prod_{d|n} \phi_d(x)$ so, by Mobius inversion, we have
$$
\phi_n(x) = \prod_{d|n}  (x^d-1)^{\mu(n/d)} .
$$
Homogenizing, we have $x_n=(r^n-s^n)/(r-s)=\prod_{d|n,\ d>1}
\phi_d(r,s)$ where $\phi_n(r,s):=s^{\phi(n)} \phi_n(r/s)\in
\Bbb Z[r,s]$. Indeed for any Lucas sequence $\{ x_n\} $ the numbers $\phi_n$, defined by
$$
\phi_n:=\prod_{d|n}  x_d^{\mu(n/d)},
$$
are integers.  Most importantly, this definition yields that $p$
is a characteristic prime factor of $\phi_n$ if and only if $p$ is a
characteristic prime factor of $x_n$; moreover $p$ divides both
$\phi_n$ and $x_n$ to the same power. Therefore $y_n$ divides
$\phi_n$, which divides $x_n$. In fact $y_n$ and $\phi_n$ are very
close to each other multiplicatively (as we show in Corollaries 3
and 4 below): either $\phi_n=y_n$, or $\phi_n=py_n$ where $p$ is
some prime dividing $n$, in this case, $n=p^em$ where $p$ is a characteristic
prime factor of $\phi_m$. So if we can show that

(i) \ $\phi_n$ is not a square, and

(ii) \ $p\phi_n$ is not a square when $n$ is of the form $n=p^em$  where $p$ is an odd prime, $e\geq 0,\  m>1$ and $m$ divides $p-1,\ p $ or $p+1$

\noindent then we can deduce that $y_n$ is not a square.
To prove this we modify the approach of Terjanian described above: We will show that there exist  integers $k$ and $\ell$  for which $$\left(  \frac{x_k}{\phi_n} \right)=\left(  \frac{x_\ell}{p\phi_n} \right)=-1,$$ where $\left(  \frac{.}{.} \right)$ is the Jacobi symbol.

Our first step then is to evaluate the Jacobi symbol $(x_k/x_m)$
for all positive integers $m$ and $k$. In fact this equals 0 if and only if $(k,m)>1$.  Otherwise, we will show that for any coprime positive integers $k$ and $m>2$ we have
$$
\left(  \frac{x_k}{x_m} \right) = (-1)^{[2u/m]}     \tag{4}
$$
for any integer $u$ which is $\equiv 1/k \pmod m$, as discussed above. (Lenstra's observation that (4) holds when $x_m=2^m-1$, which he shared with me in an email, is really the starting point for the proofs of our main results).

From this we deduce that
$$
\left(  \frac{x_k}{\phi_m} \right) =  (-1)^{N(m,u)}   \tag{5}
$$
for all $m\geq 1$, where, for $r(m)=\prod_{p|m} p$ and the Mobius function $\mu(m)$, we have
$$
N(m,u):= \mu^2(m) +  \# \{ i: 1\leq i< 2ur(m)/m\ \text{and} \ (i,m)=1 \}   .
$$
Now if $\phi_m$ is a square then by (5), we have that $N(m,u)$ is even whenever $(u,m)=1$. In Proposition 4.1 we show that this is false unless $m=1,2$ or $6$: our proof of this elementary fact is more
complicated than one might wish.

In Lemma 5.2 we show, using (5), that if $p\phi_m$ is a square
where $m=p^en,\ n>1$ and $n$ divides $p-1, p$ or $p+1$ then
$N(m,u')-N(m,u)$ is even whenever  $u\equiv u' \pmod n$ with
$(uu',m)=1$. In Propositions 5.3 and 5.5 we show that this is
false unless $m=6$: again our proof of this elementary fact is
more complicated than one might wish.

Since $x_d\equiv 3 \pmod 4$ for all $d\geq 2$ (as may be proved by induction), and since any squarefree integer $m$ has exactly $2^\ell-1$ divisors $d>1$, where $\ell$ is the number of prime factors of $m$, therefore $\phi_m \equiv \prod_{d|m} x_d \equiv x_1 3 \equiv 3 \pmod 4$, and so cannot be a square. Hence neither
$\phi_2$ nor $\phi_6$ is a square (despite the fact that
$(x_k/\phi_6)=1$ for all $k$ coprime to 6, since  $N(6,u)$ is even whenever $(u,6)=1$).  Therefore the only possibility left is that  $3\phi_6$ is a square, as claimed.
\enddemo

\demo{Proof of Corollary 2}  If $p$ is a characteristic prime factor of $x_{n_k}$, which divides $x_{n_k}$ to an odd power then $p$ does not divide $x_{n_i}$ for any $i<k$ and so divides $\prod_{1\leq i\leq k} x_{n_i}$ to an odd power, contradicting the fact that this is a square. Therefore $n_k=2$ or $6$ by Theorem 2.
Since a similar argument may be made for any $x_{n_i}$ where $n_i$ does does not divide $n_j$ with $j>i$ we deduce,  from Theorem 1, that every $n_i$ must divide $6$.

Therefore either $k=1$ and $x_2=b$ is a square, or we can rewrite $\prod_{1\leq i\leq k} x_{n_i}$ as a product of $\prod_{1\leq j\leq \ell} \phi_{m_j}$  times a square, where
$1<m_1<\dots <m_\ell=6$ and $\{ m_1,\dots,m_{\ell-1}\} \subset
\{ 2,3\}$.  However $\phi_3$ is divisible by some characteristic odd prime factor $p$ to an odd power, which does not divide $\phi_6$
(as all $x_n,\ n\geq 1$ are odd),
and so $\phi_3$ cannot be in our product.  Now $\phi_6$ is not a square since $\phi_6=b^2+3c\equiv 3\pmod 4$. Therefore both  $\phi_2$ and  $\phi_6$ are 3 times a square, which  is equivalent to $x_3x_6$ being a square.
\enddemo

Theorem 1 follows from Theorem 2, and Corollary 1 follows from Corollary 2.

\head 2. Elementary properties of Lucas sequences\endhead

\subhead 2a. Lucas sequences in general \endsubhead
If $y_{n+2}=-by_{n+1}+cy_n$ for all $n\geq 0$ with $y_0=0,\ y_1=1$
then $y_n=(-1)^{n-1}x_n$ for all $n\geq 0$. Therefore the prime factors, and characteristic prime factors, of $x_n$ and $y_n$ are the same and divide each to the same power, and so we may assume, without loss of generality, that $b>0$.

Let $\alpha$ and $\beta$ be the roots of $T^2-bT-c$. Then
$$
x_n=\frac{\alpha^n-\beta^n}{\alpha-\beta} \ \ \text{for all} \ \ n\geq 0
$$
(as may be proved by induction). We note that
$\alpha+\beta=b$ and $\alpha\beta=-c$, so that $(\alpha,\beta)|(b,c)=1$ and thus $(\alpha,\beta)=1$. Moreover $\Delta=(\alpha-\beta)^2=b^2+4c$.

In this subsection we prove some standard facts about Lucas sequences that can be found in many places (see, e.g. [3]).

\proclaim{Lemma 1} We establish various properties of the sequence $\{ x_n\}$:

\noindent {\rm (i)}\ We have $(x_n,c)=1$ for all $n\geq 1$.

\noindent {\rm (ii)}\ We have $(x_{n}, x_{n+1})=1$ for all $n\geq 0$.

\noindent {\rm (iii)}\ We have $x_{d+j}\equiv x_{d+1} x_j \pmod {x_d}$ for all $d\geq 1$ and $j\geq 0$. Therefore if $k-\ell=jd$. then $x_k \equiv x_\ell x_{d+1}^j \pmod {x_d}$.

\noindent {\rm (iv)}\ Suppose $d$ is the minimum integer $\geq 1$ for which $x_d$ is  divisible by given integer $r$. Then $r|x_k$ if and only if $d|k$.

\noindent {\rm (v)}\ For any two positive integers $k$ and $m$ we have $(x_k,x_m)=x_{(k,m)}$.

\endproclaim

\demo{Proof}
 (i)\  If not, select $n$ minimal so that there exists a prime $p$ with
$p|(x_n,c)$. Then $bx_{n-1}=x_n-cx_{n-2}\equiv 0\pmod p$
and so $p|x_{n-1}$ since $(p,b)|(c,b)=1$, contradicting minimality.
\smallskip
(ii)\ We proceed by induction using that $(x_{n+1}, x_{n+2})| x_{n+2}-bx_{n+1}=cx_n$, and thus   divides $x_n$, since  $(x_{n+1},c)=1$ by (i). Therefore $(x_{n+1}, x_{n+2})| (x_{n}, x_{n+1})=1$.
\smallskip
(iii)\ We proceed by induction on $j$: it is trivially true for $j=0$ and $j=1$. For larger $j$ we have $x_{d+j}=bx_{d+j-1}+cx_{d+j-2} \equiv x_{d+1}(bx_{j-1}+cx_{j-2})=x_{d+1} x_j \pmod {x_d}$.
\smallskip
(iv)\  Since $(x_{d+1},x_d)=1$ we see that $(x_d,x_{d+j})=(x_d,x_j)$ by (iii). So if $j$ is the least positive residue of $k \pmod d$ we find that
$(r,x_{k})=(r,x_j)$.  Now $0\leq j\leq d-1$ and  $(r,x_j)=r$ if and only if $j=0$, and hence $d|k$, so     the result follows by the definition of $d$.
 \smallskip
(v)\  Let  $g=(k,m)$ so (iv) implies that $x_g|(x_k,x_m)=r$, say. Let $d$ be the minimum integer $\geq 1$ for which $x_d$ is  divisible by $r$. Then  $d|(k,m)=g$ by (iv), and thus
$r|x_g$ by (iv), and the result is proved.

\enddemo

\proclaim{Proposition 1}  There exists an integer $n\geq 1$ for which prime $p$ divides $x_n$ if and only if $p$ does not divide $c$. In this case  let $q=p$ if $p$ is odd, and $q=4$ if $p=2$. Select $r_p$ to be the minimal integer $\geq 1$ for which $q| x_{r_p}$.  Define $e_p\geq 1$ so that $p^{e_p}| x_{r_p}$ but $p^{e_p+1}$ does not. Then $q|x_n$ if and only if $r_p|n$, in which case, writing $n=r_pp^km$ where $p\not|m$ for some integer $k$, we have that $p^{e_p+k}$ divides $x_n$ but $p^{e_p+k+1}$ does not. Finally, if $p$ is an odd prime for which $p|\Delta$ then $p|x_p$, and $p^2\nmid x_p$ if $p>3$.\endproclaim

\demo{Proof} Since  $p|x_n$ for some $n\geq 1$ we have
$(p,\alpha\beta)|(x_n,c)=1$ by Lemma 1(i) so that $p$ is coprime to both $\alpha$ and $\beta$. On the other hand if $(p,\alpha\beta)=1$ then $\alpha,\beta$ are in the group of units mod $p$, and therefore there exists an integer $n$ for which
$\alpha^n\equiv 1\equiv\beta^n \pmod p$ so that $p|\alpha^n-\beta^n$. Hence $p|x_n$ if $(p,\alpha-\beta)=1$. Now
$(p,\alpha-\beta)>1$ if and only if $p|\Delta$. In this case one easily shows, by induction, that $x_n\equiv n(b/2)^{n-1} \pmod p$ if $p>2$, and hence $p|x_p$. Finally $2|\Delta$ if and only if $2|b$, whence $c$ is odd (as $(b,c)=1$) and so $x_n\equiv n \pmod 2$; in particular $2|x_2$.

Let us write $\beta^d=\alpha^d+(\beta^d-\alpha^d)$, so that
$$
\beta^{kd}=(\alpha^d+(\beta^d-\alpha^d))^k = \alpha^{kd} + k\alpha^{(k-1)d}(\beta^d-\alpha^d) + \binom k2 \alpha^{(k-2)d}(\beta^d-\alpha^d)^2 + \ldots ,
$$
and therefore, since $x_d$ divides $x_{kd}$,
$$
x_{kd} / x_d \ \equiv \   k\alpha^{(k-1)d}  + \binom k2 \alpha^{(k-2)d}(\beta-\alpha)x_d \pmod {x_d^2} .
$$
We see that if $p|x_d$ then $p| x_{kd} / x_d$ if and only if $p|k$, as $(p,\alpha)=1$ (since $\alpha|c$ and $(p,c)=1$ by Lemma 1(i)).  We also deduce that
$
x_{pd} / x_d \ \equiv \   p\alpha^{(p-1)d}  \pmod {p^2} ,
$
and so $p^2\nmid x_{pd} / x_d$, unless $p=2$ and $x_d\equiv 2 \pmod 4$.  The result then follows from Lemma 1(iv).

Finally, if odd prime $p|\Delta=(\alpha-\beta)^2$ then
$$
x_p=\frac{\beta^{p}- \alpha^{p}}{\beta-\alpha} = p\alpha^{p-1}  + \binom p2 \alpha^{p-2}(\beta-\alpha) + \ldots \equiv 0 \pmod p.
$$
Therefore $n_p|p$ by Lemma 1(iv) and $n_p\ne 1$ (as $x_1=1$), and so $n_p=p$. Adding the two such identities with the roles of $\alpha$ and $\beta$ exchanged, yields
$$
\frac{2x_p}p =  \sum\Sb 1\leq j\leq p \\ j \ \text{odd} \endSb  \frac 1p \binom pj \Delta^{\frac{j-1}2}
 (\alpha^{p-j}+\beta^{p-j}) - \sum\Sb 1\leq j\leq p \\ j \ \text{even} \endSb  \frac 1p \binom pj \Delta^{\frac{j}2}
 x_{p-j} .
$$
This is $\equiv   \alpha^{p-1}+\beta^{p-1} \pmod p$ plus
$ \frac 23  \Delta$ if $p=3$. Now if $p>3$ the first term $=x_{2p-2}/x_{p-1}$ and so is not divisible by $p$. One can verify that $9|x_3$ if and only if $9|b^2+c$.

\enddemo

\proclaim{Corollary 3} Each $\phi_n$ is an integer.
When $p$ is a characteristic prime factor of $\phi_n$ define $n_p=n$. Then $p$ divides both $x_{n_p}$ and $\phi_{n_p}$ to the same power.  Otherwise if prime $p|\phi_n$ where $n\ne n_p$ then $n/n_p$ is a power of $p$, and $p^2\not|\phi_n$ with one possible exception: if $p=2$ with $b$  odd and $c\equiv 1 \pmod 4$ then $n_2=3$ and $2^2|\phi_6$. If $p$ is an odd prime for which $p^2|\Delta$ then $p|\phi_p$ but $p^2\nmid \phi_p$.
\endproclaim

\demo{Proof} Note first that $n_p=r_p$ when $p\ne 2$. We use the formula $\phi_n=\prod_{d|n} x_d^{\mu(n/d)}$.
If $n_p=n$ then $x_n$ is the only term on the left that is divisible by $p$, and so $p$ divides both $x_{n_p}$ and $\phi_{n_p}$ to the same power. To determine the power of $p$ dividing $\phi_n$ we will determine the power of $p$ dividing each $x_d$. To do this we begin by studying those $d$ for which $q$ divides $x_d$ (in the notation of Proposition 1), and then we return, at the end, to those $x_d$ divisible by $2$ but not $4$:

By Proposition 1, $q$ divides $x_d$ if and only if  $d=r_pp^\ell q$ with $0\leq \ell\leq k$ and $q|m$, and so the power of $p$ dividing these terms in our product is
$$
 \sum\Sb 0\leq \ell \leq k  \endSb  \mu( p^{k-\ell})   (e_p+\ell) \sum\Sb q|m \endSb  \mu(m/q) =\cases
1 & \text{if} \ m=1 \ \text{and} \  k\geq 1 \\
0 &  \text{if} \ m\geq 2 \\
e_p & \text{if} \ m=1 \ \text{and} \  k=0\ (\text{i.e. } n=r_p) .
\endcases
$$
Hence if $p$ is odd, or $p=2$ with $n_2=r_2$,  then $p|\phi_n$ with  $n>n_p$ if and only if $n/n_p$ is a power of $p$, and then $p^2\nmid \phi_n$.

Other $x_d$ divisible by $p$ occur only in the case that $p=2$ and $r_2=2n_2$,  and these are  the terms $x_d$ in the product
for which $n_2|d$ but $r_2$ does not. Such $x_d$ are divisible by $2$  but not $4$. Hence  the total power of $2$ dividing the product of these terms is
$$
\sum\Sb d|n \\ n_2|d,\ 2n_2\nmid d \endSb \mu(n/d)=\cases
1 & \text{if} \ n=n_2 \\
-1 &  \text{if} \ n= 2n_2 \\
0 & \text{otherwise}
\endcases
$$
We deduce that $2|\phi_n$ with  $n>n_2$ if and only if $n/n_2$ is a power of $2$. Moreover $4\nmid \phi_n$, except in the special case that $n=r_2=2n_2$ and $e_2\geq 3$. We now study this special case:\  We must have $c$ odd, else $c$ is even, so that $b$ is odd, and $x_n$ is odd for all $n\geq 1$. We must also have $b$ odd, else $x_n\equiv n \pmod 2$, so $n_2=2$, that is $x_2=b$ is divisible by 2 but not 4. But then $r_2=4$ and so
$\phi_4=b^2+2c\equiv 2 \pmod 4$,  a contradiction. In this case $n_2=3$ and we want $r_2=6$.  But then
$\phi_3=b^2+c\equiv 2 \pmod 4$, so that $c\equiv 2-b^2\equiv 1 \pmod 4$, and $\phi_6=b^2+3c\equiv 1+3\equiv 0 \pmod 4$.

The last statement follows from the last of Proposition 1 since $\phi_p=x_p$ (and working through the possibilities when $p=3$).

 \enddemo

Since $\phi_n$ is usually significantly smaller than $x_n$ and since we have a very precise description of the non-characteristic prime factors of $\phi_n$,  it is easier to study characteristic prime factors of $x_n$ by studying the factors of $\phi_n$

\proclaim{Lemma 3} Suppose that $p$ is a prime that does not divide $c$  (so that $n_p$ exists). Then $n_p\leq p+1$. Moreover if  $p>2$ then $n_p$ divides $p-(\Delta/p)$. \endproclaim

\demo{Proof} Proposition 1 implies this when $p|\Delta$.
We have $\alpha=(b+\sqrt{\Delta})/2$ and $\beta=(b-\sqrt{\Delta})/2$, which implies that
$$
\alpha^p \equiv \frac{b^p+\sqrt{\Delta}^p}{2^p} \equiv
\frac{b+\Delta^{(p-1)/2}\sqrt{\Delta}}{2} \equiv \frac{b+(\Delta/p)\sqrt{\Delta}}{2}
\pmod p ,
$$
and analogously $\beta^p\equiv (b-(\Delta/p)\sqrt{\Delta})/2$. Hence
if $(\Delta/p)=-1$ then $\alpha^p\equiv \beta \pmod p$ and $\beta^p\equiv \alpha \pmod p$, so that
$\alpha^{p+1} = \alpha \alpha^{p} \equiv
\alpha\beta=-c \pmod p$ and similarly $\beta^{p+1}\equiv -c \pmod p$. Now $(\alpha-\beta,p)|(\Delta,p)=1$  and therefore $p|x_{p+1}$. If $(\Delta/p)=1$ then $\alpha^{p-1} = \alpha^{-1} \alpha^{p} \equiv
\alpha^{-1} \alpha=1 \pmod p$ and similarly $\beta^{p-1}\equiv 1 \pmod p$, so that $p|x_{p-1}$.

In the special case that $p=2$ we have $c$ is odd. We see easily that if $b$ is even (and so $2|\Delta$) then $n_2=2$. If $b$ is odd then $n_2=3$ and $b^2+4c\equiv 1+4=5 \pmod 8$. Therefore
$n_2$ divides $2-(\Delta/2)$, with the latter properly interpreted.
 \enddemo

\proclaim{Corollary 4} Each $\phi_n$ has at most one non-characteristic prime factor, except $\phi_6$ is divisible by 6 if $b\equiv 3 \pmod 6$ and $c\equiv 1 \pmod 2$, and $\phi_{12}$ is divisible by 6 if $b\equiv \pm 1 \pmod 6$ and $c\equiv 1 \pmod 6$.
\endproclaim

\demo{Proof} Suppose  $\phi_n$ has two non-characteristic prime factors $p<q$. By Corollary 3 we have that $q|n_p$ and so $q\leq n_p\leq p+1$ by Lemma 3. Therefore $p=2$ and $q=3$, in which case $n_2=3$, so that
$n=2^e3$ for some $e\geq 1$, and this equals $3^fn_3$ for some $f\geq 1$ by Corollary 3. Thus $f=1$ and $n_3=2$ or $4$.
The result follows by working through the possibilities mod 2 and mod 3.
\enddemo

\proclaim{Corollary 5} Suppose that $x_n$ does not contain a characteristic prime factor to an odd power and $n\ne 6$ or $12$. Then either $\phi_n=\square$ (where $\square$ represents the  square of an integer), or $\phi_n=p\square$ where $p$ is a prime for which $p^e|n$ with $e\geq 1$ and $n/p^e\leq p+1$.
\endproclaim

\demo{Proof} Follows from Corollaries 3 and 4 and Lemma 3.
\enddemo

\proclaim{Lemma 4} Suppose that the odd prime $p$ divides $\Delta$. Then $x_n\equiv n(b/2)^{n-1} \pmod p$ for all $n\geq 0$.
\endproclaim

\demo{Proof} This follows by induction on $n$: it is trivially true for $n=0,1$, and then
$$\align
x_n&=bx_{n-1}+cx_{n-2} \equiv b (n-1)(b/2)^{n-2}+c(n-2)(b/2)^{n-3} \\
&\equiv
2 (n-1)(b/2)^{n-1}-(n-2)(b/2)^{n-1} = n(b/2)^{n-1} \pmod p,
\endalign
$$
since $\Delta=b^2+4c\equiv 0 \pmod p$, so that $c\equiv -(b/2)^2 \pmod p$.
\enddemo

\subhead 2b. Lucas sequences with $b, \Delta>0,\ (c/b)=1$ and $b\equiv 3 \pmod 4, \ c\equiv 2\pmod 4$ \endsubhead
As $b, \Delta>0$ this implies that $x_n>0$ for all $n\geq 1$
since $\alpha>|\beta|$.

We also have  $x_n\equiv 3 \pmod 4$ for all $n\geq 2$, by induction.  In fact  $x_{n+2}\equiv x_n \pmod 8$ for all $n\geq 3$, which we can prove by induction: We have $$x_5=b^4+3cb^2+c^2\equiv 1+3c+4\equiv 1+c \equiv b^2+c=x_3 \pmod 8,$$  and  $$x_6=b(b^4+4cb^2+3c^2)\equiv b(1+0+4)= b(1+4)\equiv  b(b^2+2c)=x_4 \pmod 8.$$ For larger $n$, we then have $x_{n+2}=bx_{n+1}+cx_n\equiv bx_{n-1}+cx_{n-2}=x_n \pmod 8$ by the induction hypothesis.

We also note that $x_{n+2}\equiv bx_{n+1} \pmod {c}$ for all $n\geq 0$, and so
$x_n\equiv b^{n-1} \pmod c$ for all $n\geq 1$. We deduce from this and the previous paragraph that
$x_{n+2} \equiv b^2 x_n \pmod {4c}$ for all $n\geq 3$.

\proclaim{Proposition 2}  We have $(x_{d+1}/x_d)=1$ for all $d\geq 1$.
\endproclaim

\demo{Proof} For $d=1$ this follows as $x_1=1$; for $d=2$ we have
$(x_3/x_2)=( (b^2+c)/b)=(c/b)=1$. The result then follows from proving that  $\theta_d:=(x_{d+1}/x_d)(x_d/x_{d-1})=1$ for all $d\geq 3$.
Since $x_{d+1}\equiv cx_{d-1} \pmod {x_d}$ and as $x_d\equiv x_{d-1}\equiv 3\pmod 4$ for $d\geq 3$, we have
$\theta_d=(cx_{d-1}/x_d)(x_d/x_{d-1}) =-(c/x_d)=(-c/x_d)$. We will prove that this equals $1$ by induction on $d\geq 3$. So write $-c=\delta C$ where $C=|c/2|$. Then  note that
$$
\theta_3=\left( \frac{-c}{b^2+c} \right)=\left( \frac{\delta}{b^2+c} \right) \left( \frac{C}{b^2+c} \right)= \left( \frac{\delta}{b^2+c} \right) \left( \frac{-1}{C} \right)  \left( \frac{b^2+c}{C} \right)=\left( \frac{\delta}{b^2-\delta C} \right)  \left( \frac{-1}{C} \right)
$$
which is shown to be $1$, by running through the possibilities $\delta=\pm 2$ and $C\equiv \pm 1 \pmod 4$. Also, as $(-c/b)=-1$,
$$
\theta_4=\left( \frac{-c}{b(b^2+2c)} \right)=-\left( \frac{\delta}{b^2+2c} \right) \left( \frac{C}{b^2+2c} \right)=-(-1)  \left( \frac{b^2+2c}{C} \right)  =1
$$
since $\delta=\pm 2$ and $b^2+2c\equiv 5 \pmod 8$. Now for the induction hypothesis, for $d\geq 5$:\
The value of $\theta_d=(-c/x_d)$ depends only on the square class of $x_d \pmod {4c}$, and we saw in the paragraph above that this is the same square class as $x_{d-2} \pmod {4c}$ for $d\geq 5$. Hence $\theta_{d}=1$ for all $d\geq 3$, and the result follows.

\enddemo

\head 3. Evaluation of Jacobi symbols, when $b, \Delta>0,\ b\equiv 3 \pmod 4, \ c\equiv 2\pmod 4$ and $(c/b)=1$. \endhead

\subhead 3a. The reciprocity law \endsubhead
Suppose that $k$ and $m>1$ are coprime positive integers.
Let $u_{k,m}$ be the
least residue, in absolute value, of $1/k \pmod m$
(that is $u\equiv k \pmod m$ with $-m/2< u\leq m/2$).

\proclaim{Lemma 5} If $m,k\geq 2$ with $(m,k)=1$ then
$ku_{k,m} +mu_{m,k}=1$.
\endproclaim

\demo{Proof} Now
$v:=(1-ku_{k,m})/m$ is an integer $\equiv 1/m \pmod k$ with
$-k/2 +1/m \leq v < k/2 +1/m$. This implies that $-k/2<v\leq k/2$, and so $v=u_{m,k}$.
\enddemo

\proclaim{Theorem 4} If $k\geq 1$ and $m>1$ are coprime positive integers then the value of the Jacobi symbol $(x_k/x_m)$ equals the sign of $u_{k,m}$.
\endproclaim

\demo{Proof} By induction on $k+2m\geq 5$.  Note that when $k=1$ we have $u=1$ and the result follows as $(x_1/x_m)=(1/x_m)=1$.  For larger $k$, we have two cases. If $k>m$ then let $\ell$ be the least positive residue of $k \pmod m$, say $k-\ell=jm$. By Lemma 1(iii) we have
$(x_k/x_m)=(x_\ell/x_m)(x_{m+1}/x_m)^j=(x_\ell/x_m)$ by Proposition 2.  Moreover
$u_{l,m}=u_{k,m}$ by definition so that the result follows from the induction hypothesis. If $2\leq k<m$ then
$(x_k/x_m)=-(x_m/x_k)$ since $x_m\equiv x_k\equiv 3\pmod 4$.
Moreover $u_{k,m}$ and $u_{m,k}$ must have opposite signs, else
$1=k|u_{k,m}| +m|u_{m,k}|\geq 1+1$ by Lemma 5 which is impossible. The result follows  from the induction hypothesis.
\enddemo

Define $(t)_m$ to be the least (positive) residue of $t \pmod m$,
so that $(t)_m=t-m[t/m]$. Note that $0\leq (t)_m< m/2$
if and only if $[(t)_m/(m/2)]=0$. Also that
$[(t)_m/(m/2)]=[2t/m]-2[t/m] \equiv [2t/m] \pmod 2$
Now, if $m\geq 3$ and $(t,m)=1$ then $(t)_m$
is not equal to $0$ or $m/2$; therefore
if $u$ is {\sl any} integer $\equiv 1/k \pmod m$ then  the sign of $u_{k,m}$ is given by $(-1)^{ [2u/m] }$. We deduce the following from this and Theorem 4:

\proclaim{Corollary 6}  Suppose that $k$ and $m\ne 2$ are coprime positive integers. If $u$ is any integer $\equiv 1/k \pmod m$ then
$$
\left(  \frac{x_k}{x_m} \right) = (-1)^{[2u/m]} .   \tag{4}
$$
\endproclaim

Note that if $k$ is odd then $ \left(  \frac{x_k}{x_2} \right) =1$, whereas (4) would always give $-1$.

\remark{Remark} In email correspondence with Ilan Vardi we understood how (4) can be deduced directly from (3) and known facts about continued fractions. Write $p_n/q_n= [a_0,a_1,\dots, a_n]$ for each $n$, and recall that
$$
\left( \matrix p_n&p_{n-1}\cr q_n&q_{n-1} \cr  \endmatrix \right) =
\left( \matrix a_0&1\cr 1&0 \cr  \endmatrix \right) \left( \matrix a_1&1\cr 1&0 \cr  \endmatrix \right) \dots
\left( \matrix a_{n}&1\cr 1&0 \cr  \endmatrix \right)
$$
as may easily be established by induction on $n\geq 1$. By taking determinants we see that $p_nq_{n-1}=  p_{n-1}q_n+(-1)^{n+1}\equiv (-1)^{n+1} \pmod {q_n}$.
Taking $p_{n}/q_{n}=k/m$ with $n=\Lambda(k/m)-1$ and
$u$ to be the least positive residue of $1/k \pmod m$ we see that $q_{n-1}\equiv (-1)^{n+1}u \pmod {m}$ and
$q_{n-1}<q_n=m$, so
$q_{n-1}=u$ if $n$ is odd, $q_{n-1}=m-u$ if $n$ is even.
Now $m=q_n=a_nq_{n-1}+q_{n-2}\geq 2q_{n-1}+1$, and so
$q_{n-1}<m/2$. Therefore if $u<m/2$ then $q_{n-1}=u$, so $n$ is odd and the values given in (4) and (3) are equal. A similar argument works if $u>m/2$. Hence we have that
$$
\Lambda(k/m) \equiv [2u/m] \pmod 2 \ \ \text{where} \ \ uk\equiv 1 \pmod m \tag{6}
$$
for all coprime, positive integers $k$ and $m$.
\endremark

\subhead 3b.  The  characteristic part \endsubhead
If $(m,k)=1$ and $u\equiv 1/k \pmod m$ then
$$
\left(  \frac{x_k}{\phi_m} \right)   = \prod_{d|m} \left(  \frac{x_k}{x_d} \right)^{\mu(m/d)} =  (-1)^{E(m,u)} \tag{7}
$$
by (4) since $(x_k/x_d)=1$ if $d=1$ or $2$, where
$$
\align
E(m,u) &\equiv  \sum\Sb d|m\\ d\geq 3 \endSb \mu\left( \frac md \right) \left[ \frac{2u}{d} \right] = \sum\Sb d|m\\ d\geq 3 \endSb \mu\left( \frac md \right) \sum\Sb {1\leq j\leq 2u-1} \\ d|j \endSb 1 \\
&\equiv
\sum\Sb {1\leq j\leq 2u-1}  \endSb \sum\Sb d|(m,j) \endSb \mu\left( \frac md \right) + \mu(m)(2u-1) +E_2 \pmod 2
\endalign
$$
where $E_2$, the
contribution when $d=2$,  occurs only when $m$ is even, and is then  equal to $\mu(m/2)(u-1)$, and we can miss the $j=2u$ term since if $d|2u$ then $d|(2u,m)=(2,m)|2$. However $u$ is then odd since $(u,m)=1$ and so  $E_2\equiv \mu(m/2)(u-1)\equiv 0 \pmod 2$.

Now let $r(n)=\prod_{p|n} p$ for any integer $n$. We see that
$\mu(m/d)=0$ unless $m/d$ divides $r(m)$, that is $d$ is divisible by $m/r(m)$, in which case $j$ must be also.
Write $j=i(m/r(m))$, and each $d=D(m/r(m))$ and so
$$
E(m,u) \equiv \mu(m) +
\sum\Sb {1\leq i< \frac{2ur(m)}m}  \endSb \ \sum\Sb D|(r(m),i) \endSb \mu(r(m)/D ) \equiv  \mu(m) +
\sum\Sb {1\leq i< 2ur(m)/m} \\ (i,m)=1  \endSb \ 1
\pmod 2     ,
$$
which is $N(m,u)$, and so we obtain (5).

\head 4.  The tools needed to show that $\phi_m\ne \square$. \endhead

\proclaim{Proposition 4.1} If $m\ne 1,2,6$ then $N(m,u')-N(m,u)$ is odd for some $u, u'$ with $(uu',m)=1$.
\endproclaim

\demo{Proof}
If $m$ is squarefree then
$N(m,u')-N(m,u)= \# \{ i: 2u\leq i< 2u'\ \text{and} \ (i,m)=1 \}$. So, if  $m$ is odd and $>1$
let $u=(m-1)/2$ and $u'=u+1$.
If $m$ is even then there exists
a prime $q|m$ with $q\geq 5$ (as $m\ne 2$ or $6$), so we can write $m=qs$ where $q\nmid s>1$: Then select $u\equiv -1 \pmod {s}$ and $u\equiv -3/2 \pmod{q}$ with $u'=u+2$.
\smallskip

For $m$ not squarefree let $m_2$ be the largest powerful number dividing $m$ and $m=m_1m_2$ so that $m_1$ is squarefree,
$(m_1,m_2)=1$, and $r(m_2)^2|m_2$. Note that $m/r(m)=m_2/r(m_2)$.

\smallskip

When $m_2=4$ then $N(m,u)= \# \{ i: 1\leq i< u,\ (i,m)=1\}$, so if $u$ is the smallest integer $>1$ that is coprime with $m$ then $N(m,u)-N(m,1)=1$.

\smallskip

So we may assume that $m_2>4$, in particular that $2r(m)/m\leq 2/3$. Consider

\noindent
$N(m, \frac m{r(m)}(\ell+1)+1 )- N(m, \frac m{r(m)}\ell+1)=
\#\{ i: 2\ell+1\leq i\leq 2\ell+2:\ (i,m)=1\}$.

\noindent Select $\ell\equiv -1 \pmod {m_2}$ so that $(2\ell+2,m)\geq m_2$. Then we need to select $\ell \pmod p$
for each prime $p$ dividing $m_1$ so that each of
$\frac m{r(m)}(\ell+1)+1,\  \frac m{r(m)}\ell+1$ and $2\ell+1$ are
coprime to $p$. Since there are just three linear forms, such congruence classes exist modulo primes $p>3$ by the pigeonhole principle; and also for $p=3$ as may be verified by a case-by-case analysis. Thus the result follows when $m_1$ is odd.

So we may assume that $m_1$ is even and now consider

\noindent
$N(m, \frac {2m}{r(m)}(\ell+1)+1 )- N(m, \frac {2m}{r(m)}\ell+1)=
\#\{ i: 4\ell+1\leq i\leq 4\ell+4:\ (i,m)=1\}$.

\noindent Select $\ell\equiv -3/4 \pmod {m_2}$ so that $(4\ell+3,m)\geq m_2$. We can again select
$\ell \pmod p$ for each prime $p>3$ dividing $m_1$ so that each of $\frac {2m}{r(m)}(\ell+1)+1,\  \frac {2m}{r(m)}\ell+1, \ 4\ell+1$ are coprime to $p$ by the pigeonhole principle, and therefore the result follows if $3$ does not divide $m_1$.

So we may assume that $6|m_1$. Select integer $\ell$ so that  $\ell\equiv 1 \pmod {m_2}, \ \ell\equiv -m/r(m) \pmod 4$ and,
for each prime $p$ dividing $m_1/2$, $p$ does not
divide $\ell, \  \frac {m}{r(m)} \ell-1$ or $\frac {m}{r(m)} \ell+3$. Therefore, since $3r(m)/m\leq 3/5$, we have

\noindent
$N(m, \frac 12(\frac {m}{r(m)} \ell+3))- N(m, \frac 12(\frac {m}{r(m)} \ell-1))=
\#\{ i:  \ell \leq i< \ell+1:\ (i,m)=1\}=1$.
\enddemo

\head 5.  The tools needed to show that $\phi_m\ne p\ \square$. \endhead

\proclaim{Lemma 5.1} Suppose that $\phi_m=p\square$, where $p$ is an odd prime,  $m=p^en,\ 1<n\leq p+1$ and $p|\phi_n$.
If $k\equiv k' \pmod {2n}$ with $(kk',m)=1$ then $(x_k/\phi_m)=(x_{k'}/\phi_m)$. Moreover if $c\equiv 2 \pmod 4$ then $(\phi_m/x_k)=(\phi_m/x_{k'})$.
\endproclaim

\demo{Proof} Writing $k'=k+2nj$ we have
$x_{k'}\equiv x_k x_{n+1}^{2j} \pmod {x_n}$, by Lemma 1(iii); and so $(x_{k}/p)=(x_{k'}/p)$ since $p|x_n$. Therefore since $\phi_m=p\square$ we have
$(x_k/\phi_m)=(x_k/p)=(x_{k'}/p)=(x_{k'}/\phi_m)$.

If $c\equiv 2 \pmod 4$ and $k\equiv k' \pmod 2$ then $x_k\equiv x_{k'} \pmod 4$, which implies that $\left(  {p}/{x_k} \right) \left(  {p}/{x_{k'}} \right) = \left(  {x_k}/{p} \right) \left(  {x_{k'}}/{p} \right)$, and the result follows from the first part.
\enddemo

\proclaim{Lemma 5.2} Assume that $b, \Delta>0,\ b\equiv 3 \pmod 4, \ c\equiv 2\pmod 4$ and $(c/b)=1$. Suppose that $\phi_m=p\square$, where $p$ is an odd prime,  $m=p^en,\ 1<n\leq p+1$ and $p|\phi_n$.
If $u\equiv u' \pmod n$ with $(uu',m)=1$ then
$N(m,u')-N(m,u)$ is even. If $e=1$ and $n\ne p$ then this implies that $N(n,u'/p)-N(n,u/p)$ is even.
\endproclaim

\demo{Proof} Let $k, k^*$ be integers for which $k\equiv 1/u \pmod m$ and $k^*\equiv 1/u' \pmod m$. Evidently $k\equiv 1/u \equiv 1/u' \equiv k^*  \pmod n$.
If $k\equiv k^* \pmod {2n}$ then let
$k'=k^*$, otherwise take $k'=k^*+m$, so $k'\equiv k  \pmod {2n}$
(since $m/n=p^e$ is odd). Applying the first part of Lemma 5.1, the first result follows from (5).

If $e=1$ then $m=pn$ so that $r(m)/m=r(n)/n$.  Therefore $N(m,u')-N(m,u)$ equals, for $U= 2ur(n)/n$ and $U'= 2u'r(n)/n$,
$$
\sum\Sb U \leq i< U' \\ (i,r(n)p)=1 \endSb 1
= \sum\Sb U \leq i< U' \\ (i,r(n))=1 \endSb 1 - \sum\Sb U \leq i< U' \\ (i,r(n))=1,\ p|i \endSb 1 \equiv \sum\Sb U/p \leq j< U'/p \\ (j,r(n))=1  \endSb 1   \pmod 2.
$$
since $U'\equiv U \pmod {2r(n)}$ (as $u\equiv u' \pmod n$),
so that the  first term counts each residue class coprime with $r(n)$ an even number of times, and by writing $i=jp$ in the second sum.  The result follows.
\enddemo

\proclaim{Proposition 5.3} Suppose $n\geq 2$ and  $n$ divides $p-1$ or $p+1$ for some odd prime $p$. Let  $m=p^en$ for some $e\geq 1$. There exists an integer $u$ such that $(u(u+n),m)=1$ for which
$N(m,u+n)-N(m,u)=1$ if $e\geq 2$, for which $N(n,(u+n)/p)-N(n,u/p)=1$ if $e=1$, except when $p=3, n=2$. In that case we have $N\left(2\cdot 3^e, \frac{3^{e-1}+4+3(-1)^e}2 \right)-N\left(2\cdot 3^e,1\right)=1$, for $e\geq 2$.
\endproclaim

\proclaim{Lemma 5.4} If $n\geq 3$ and odd prime $p=n-1$ or $p\geq n+1$ (except for the cases $n=3$ or $6$ with $p=5$; and $n=4, p=3$) then in any non-closed interval of length $n$, containing exactly $n$ integers, there exists an integer $u$ for which $u$ and $u+n$ are both prime to $np$.
\endproclaim

\demo{Proof} Since $p\geq n-1$ there are no more than three integers, in our two consecutive intervals of length $n$, that are divisible by $p$ so the result follows when $\phi(n)\geq 4$. Otherwise $n=3,4$ or $6$, and if the reduced residues are $1<a<b<n$ then $p$ divides $b-a,\ (n+b)-a,\ (n+a)-b$ or $(2n+a)-b$. Therefore
$p|4, 10, 2$ or 8 for $n=6$;\ $p|2$ or 6 for $n=4$; $p|1, 4, 2$ or $5$ for $n=3$. The result follows.
\enddemo

\demo{Proof of Proposition 5.3}  Let $f:=\max\{ 1,e-1\}$. The result holds for
$$
(m,u)=\left(3\cdot 5^e, \frac{5^f-3}2 \right), \
\left(6\cdot 5^e,\frac{5^f-3}2 \right), \
(4\cdot 3^e,3^f-2) , \ \left(2\cdot p^e, \frac{p^f-j}2 \right)
$$
for each $e\geq 1$ and, in the last case, any prime $p>3$, where $j$ is either 1 or 3, chosen so that $u$ is odd.

Otherwise we can assume the hypotheses of Lemma 5.4.
Now suppose that $e\geq 2$. Given an integer $\ell$ we can select $u$ in the range  $\ell \frac m{2r(m)} -n < u\leq \ell \frac m{2r(m)}$ (which is an interval of length $n$) such that $u$ and $u':=u+n$ are both prime to $np$, by
Lemma 5.4.  Therefore
$N(m,u')-N(m,u)$ counts the number of integers, coprime with $m$, in an interval of length $\lambda:=2nr(m)/m=2r(n)/p^{e-1}$.
Note that $\lambda\leq  2n/p \leq 2(p+1)/p<3$ so our interval contains no more than $[\lambda]+1\leq 3$ integers, one of which is $\ell$.   If $\lambda<2$ we select $\ell\equiv 1 \pmod p$ and $\ell\equiv -1\pmod n$ so that $N(m,u')-N(m,u)=1$. Otherwise $\lambda\geq 2$ so that $n\geq r(n)\geq p^{e-1}\geq p$, and thus $n=p+1,\ e=2$ and $r(n)=n$, that is
$n$ is squarefree, and $2|(p+1)|n$. So select $\ell$ to be an
odd integer for which
$\ell\equiv 2 \pmod p$ and  $\ell\equiv -2\pmod {n/2}$
so that $\ell\pm 2, \ell\pm 1$ all have common factors with $m$,
and therefore $N(m,u')-N(m,u)=1$.

For $e=1$ and given integer $\ell$ we now select $u$ in the range  $\ell \frac {pn}{2r(n)} -n < u\leq \ell \frac {pn}{2r(n)}$, and $N(n,u'/p)-N(n,u/p)$ counts the number of integers, coprime with $n$, in an interval of length $\lambda:=2r(n)/p$. If $\lambda<1$ we select $\ell$ so that it is coprime with $n$ then we have that $N(n,u'/p)-N(n,u/p)=1$  is odd. If $\lambda\geq 1$ we have $r(n)\geq p/2$, and we know that $r(n)|n|p\pm 1$, so that $r(n)$ and $n$ equal  $\frac {p+1}2,\ p-1$ or $p+1$. If $n=r(n)=p-1$ then $n$ is squarefree and divisible by 2, and $[\lambda]=1$; so we select $\ell\equiv 1 \pmod 2$ and $\ell\equiv -1\pmod {n/2}$ so that $N(n,u'/p)-N(n,u/p)=1$. In all the remaining cases, one may check that $N(n,(n+1)/p)-N(n,1/p)=1$.
\enddemo

\proclaim{Proposition 5.5} If $m=p^{e+1}$ where  $p$ is an odd prime then
$N\left(m,\frac{p^e+1}2\right)-N(m,1)=1$.
\endproclaim

\head 6.  Other Lucas sequences \endhead

\proclaim{Proposition 6.1} Assume that $\Delta$ and $b$ are
positive with $(b,c)=1$. For $n>1$ odd with $(m,n)=1$ we have the
following:
$$
\left(  \frac{x_m}{x_n} \right) = \cases \left( \frac{c}{b}
\right)^{(m-1)(n-1)/2} & \ \text{if} \ 4|c \cr (-1)^{\Lambda(m/n)+
\left( \frac{b+1}2 \right) (m-1)} \left( \frac{c}{b}
\right)^{(m-1)(n-1)/2} & \ \text{if} \ c\equiv 2 \pmod 4 \cr
\left(  \frac{m}{n} \right)^{\frac{c-1}2} \left(  \frac{2}{n}
\right)^{(m-1)(\frac{b+c-1}2)}\left( \frac{b}{c}
\right)^{(m-1)(n-1)/2} & \ \text{if} \ 2|b \cr
\endcases
$$
\endproclaim

\demo{Proof}  For $m$ odd this is the result of Rotkiewicz [5],
discussed in section 1c. Note
that if $c$ is even then $b$ is odd and $x_n$ is odd for all
$n\geq 1$; and if $b$ is even then $c$ is odd and $x_n\equiv n
\pmod 2$ is odd for all $n\geq 1$. Thus $x_n$ is odd if and only
if $n$ is odd.

For $m$ even and $n$ odd we have that $m+n$ is odd and so
$$
\left(  \frac{x_m}{x_n} \right) = \left(  \frac{x_{m+n}}{x_n}
\right) \left(  \frac{x_{n+1}}{x_n} \right)
$$
by Lemma 1(iii); and therefore
$$
\left(  \frac{x_{n+1}}{x_n} \right) = \left(  \frac{x_2}{x_n}
\right) \left(  \frac{x_{n+2}}{x_n} \right)    ;
$$
note that $n,n+2$ are both odd, so we have yet to determine only
$(x_2/x_n)=(b/x_n)$.

Suppose that $c$ is even so that $b$ is odd. If $4|c$ then
$x_n\equiv 1 \pmod 4$ if $n$ is odd so that $(b/x_n)=(x_n/b)$. If
$c\equiv 2 \pmod 4$ and  $b\equiv 1 \pmod 4$ then
$(b/x_n)=(x_n/b)$. Now $x_n\equiv cx_{n-2} \pmod b$ and so
$x_n\equiv c^{(n-1)/2}  \pmod b$ for every odd $n$.  Therefore
$$
\left(  \frac{x_m}{x_n} \right) = \left(  \frac{x_{m+n}}{x_n}
\right) \left(  \frac{x_{n+2}}{x_n} \right)  \left(  \frac{c}{b}
\right)^{(n-1)/2} .
$$
The results follow in these cases since
$\Lambda((m+n)/n)=\Lambda(m/n)$ as $(m+n)/n=1+m/n$, and
$\Lambda((n+2)/n)=3$ as $(n+2)/n=[1,\frac {n-1}2,2]$.

If $c\equiv 2 \pmod 4$ and  $b\equiv 3 \pmod 4$ then $x_n\equiv
3\pmod 4$ for all $n\geq 2$. Therefore $(b/x_n)=-(x_n/b)$ for all
odd $n>1$, and the result follows.

Now assume that $b$ is even so that $c$ is odd. As $x_n\equiv
c^{(n-1)/2} \pmod {[b,4]}$ for each odd $n$ we have, writing
$b=2^eB$ with $B$ odd,
$$
\left(  \frac{b}{x_nc^{(n-1)/2}} \right) = \left(
\frac{2}{x_nc^{(n-1)/2}} \right)^e \left(
\frac{{x_nc^{(n-1)/2}}}{B} \right) = \left(
\frac{2}{x_nc^{(n-1)/2}} \right)^e .
$$
Now if $4|b$ then $x_n\equiv c^{(n-1)/2} \pmod 8$. Finally, if $e=1$ then $x_nc^{(n-1)/2} \equiv
1 \pmod 8$ if $n\equiv \pm 1 \pmod 8$, and $\equiv 5 \pmod 8$ if
$n\equiv \pm 3 \pmod 8$. Therefore $\left(
\frac{2}{x_nc^{(n-1)/2}} \right) =\left(  \frac{2}{n} \right)$.
The result follows.
\enddemo

\proclaim{Corollary 6.2} Suppose that $\Delta$ and $b$ are
positive, with $(b,c)=1$ and $c\equiv 2 \pmod 4$. For  $n>1$ odd, $m>1$ and $(m,n)=1$. Suppose that $mu\equiv 1 \pmod n$. If $n$ is
a power of a prime $p$ then
$$
\left(  \frac{x_m}{\phi_n} \right) =  (-1)^{N(n,u)+\mu(n)\left( \frac{b+1}2 \right) (m-1) }  \left(
\frac{c}{b} \right)^{(m-1)(p-1)/2} .
$$
If $n$ has at least two distinct prime factors then
$$
\left(  \frac{x_m}{\phi_n} \right) =  (-1)^{N(n,u)+\mu(n)\left( \frac{b+1}2 \right) (m-1) } .
$$
If $m$ is even and $>2$ then, for $nv\equiv 1 \pmod m$,
$$
\left(  \frac{\phi_m}{x_n} \right) =  (-1)^{N(m,v)+\mu(m/2)}.
$$
\endproclaim

\demo{Proof} Throughout we assume that $n>1$ is odd.
Proposition 6.1 yields that we have $\left(  \frac{x_m}{\phi_n} \right)=(-1)^A(c/b)^B$ where
$B$ equals $(m-1)/2$ times
 $$
 \sum\Sb d|n \\ d>1\endSb   \mu(n/d) (d-1) =  \sum\Sb d|n  \endSb   \mu(n/d) (d-1)
 =  \sum\Sb d|n  \endSb   \mu(n/d) d =\phi(n) \equiv \prod_{p|n} (p-1) \pmod 4
$$
which  is divisible by $4$ except if $n$ is a power of odd prime $p$, so we confirm the claimed powers of $(c/b)$. If $d|n$ then $\Lambda(m/d)\equiv [2u/d] \pmod 2$ where $um\equiv 1 \pmod n$, by (6),  and so
$$\align
A&=  \sum\Sb d|n \\ d>1\endSb   \mu(n/d) \left( \Lambda(m/d)+
\left( \frac{b+1}2 \right) (m-1) \right) \\
&\equiv
\sum\Sb d|n \\ d>1\endSb   \mu(n/d)  [2u/d] - \mu(n) \left( \frac{b+1}2 \right) (m-1) \pmod 2\\
&\equiv
N(n,u) + \mu(n) \left( \frac{b+1}2 \right) (m-1) \pmod 2
\endalign
$$
since, in section 3b, we showed that $\sum\Sb d|n\\ d\geq 3 \endSb \mu\left( \frac nd \right) \left[ \frac{2u}{d} \right] \equiv N(n,u) \pmod 2$, and here $n$ is odd (so there is no $d=2$ term).

In the third case we use the fact that if $d<n$ then the continued fraction for $d/n$ is that of $n/d$ with a $0$ on the front, and vice-versa. Hence $\Lambda(n/d)+\Lambda(d/n)\equiv 1 \pmod 2$. Hence
$$
\sum_{d|m} \mu(m/d) \Lambda(d/n) \equiv
\sum_{d|m} \mu(m/d)( \Lambda(n/d) +1)
\equiv N(m,v)+\mu(m/2) \pmod 2 .
$$
The other terms disappear since $\phi(m)$ is even.

\enddemo

\demo{Proof of Theorem 2} Our goal is to show that $y_n$ is not a square, just as we did in the proof of Theorem 2'. We begin by showing that $\phi_n$ is not a square, for $n\ne 1,2,3,6$ by using Corollary 6.2:

Suppose that $\phi_n$ is  a square so that $(x_m/\phi_n)=1$. For $n>1$ odd, we compare, in the first two identities of Corollary 6.2, the results for $m$ and $m+n$. The value of $u$ does not change and we therefore deduce that $(-1)^{\mu(n)\left( \frac{b+1}2 \right)   }  \left(
\frac{c}{b} \right)^{ (p-1)/2}=1$  and $(-1)^{\mu(n)\left( \frac{b+1}2 \right)   }  =1$, respectively. Hence those identities both become $N(n,u)\equiv 0 \pmod 2$ whenever $(u,n)=1$.
Similarly if $n>2$ is even then the third identity of Corollary 6.2 yields that $N(n,u)\equiv \mu(n/2) \pmod 2$ whenever $(u,n)=1$. These are all impossible,  by Proposition 4.1, unless  $n=1,2$ or $6$.

Next we suppose that $p\phi_n$ is a square where $n=p^em$  and $p$ is an odd characteristic prime factor of $\phi_m$, with $e\geq 0,\  m>1$ and $m$ divides $p-1,\ p $ or $p+1$. Lemma 5.1 tells us that if $k\equiv k' \pmod {2m}$ with $(kk',n)=1$ then $(x_k/\phi_n)=(x_{k'}/\phi_n)$ and
$(\phi_n/x_k)=(\phi_n/x_{k'})$. Corollary 6.2 thence implies that if $n>2$ then $N(n,u)\equiv N(n,u') \pmod 2$ where $uk\equiv u'k'\equiv 1 \pmod n$. We now proceed as in Lemma 5.2 to deduce that if $u\equiv u' \pmod m$ with $(uu',n)=1$ then $N(n,u)-N(n,u')$ is even, deduce the final part of that Lemma, and then use Proposition 5.3 to obtain the desired contradiction except for when $n=1,2$ or $6$.

We can now deduce that $y_n$ is not a square, for $n\ne 1,2,6$, from the last two paragraphs, and the result follows.
\enddemo

\demo{Proof of Theorem 3} We deduce Theorem 3 from Theorem 2 by ruling out the possibility that there exists an $n$ for which all of the characteristic prime factors $p$ of $x_n$ which divide  $x_n$ to an odd power, are not primitive prime factors of $x_n$. If this were the case then each such $p$ would be a divisor of $\Delta$, which is odd, so that $p$ is odd, and therefore $n=n_p=p$ by Lemma 3. Hence there is a unique such $p$, and we must have that $x_p=\phi_p$ is $p$ times a square. But then
$$
\left( \frac{x_m}{\phi_p}\right) = \left( \frac{x_m}{p}\right) = \left( \frac{ m (b/2)^{m-1}}{p}\right)
$$
by Lemma 4 whenever $p\nmid m$. Comparing this to the first part of Corollary 6.2 we find that
$$
\left( \frac{ m (b/2)^{m-1}}{p}\right)  =  (-1)^{N(p,u)+ \left( \frac{b+1}2 \right) (m-1) }  \left(
\frac{c}{b} \right)^{(m-1)(p-1)/2} .
$$
where $mu\equiv 1 \pmod p$. Replacing $m$ by $m+p$,  does not change $u$, so comparing the two estimates yields that $\left(  (b/2)/p \right) =  (-1)^{ \frac{b+1}2  } \left( c/b \right)^{ (p-1)/2}$ and thus the last equation becomes
$$
\left( \frac{ u  }{p}\right) =\left( \frac{ m  }{p}\right)=  (-1)^{N(p,u)}   =  (-1)^{[2u/p]}
$$
for $u\ne 1$, since $N(p,u)\equiv [2u/p] \pmod 2$ if $p\nmid u$.  Now, selecting $u=2$ we deduce that $(2/p)=1$ if $p>3$. Taking $u=\frac{p-1}2$ we obtain $(\frac{p-1}2/p) = 1$ and, taking $u= p-1 $ we obtain  $( (p-1)/p)=-1$. These three estimates imply $1\times 1=-1$, a contradiction, for all $p>3$.

\enddemo

We note that in the other cases with $bc$ even, our argument will not yield such a general result about characteristic prime factors:

\proclaim{Corollary 6.3} Suppose that $4|c$ and $b\equiv 1 \pmod 2$, with $(m,n)=1$.
Suppose that $n$ is odd:
If $n$ is a power of a prime $p$
then
$$
\left(  \frac{x_m}{\phi_n} \right) =  \left(  \frac{c}{b}
\right)^{(m-1)(p-1)/2} .
$$
Otherwise $(x_m/\phi_n)=1$ if $n$ has at least two distinct prime
factors. On the other hand if $n$ is even and $>2$ then  $(\phi_n/x_m)=1$.
\endproclaim

One can deduce that $\phi_{p^k}$ is not a square if  $4|c$ and $(c/b)=-1$ and $p\equiv 3 \pmod 4$.

\proclaim{Corollary 6.4} Suppose that $b$ is even and $c$ is odd, with $(m,n)=1$.
Suppose that $n$ is odd: If $n$ is a power of a prime $p$ then
$$
\left(  \frac{x_m}{\phi_n} \right) = \left(  \frac{m}{p}
\right)^{\frac{c-1}2} \left(  \frac{2}{p}
\right)^{(m-1)(\frac{b+c-1}2)}\left( \frac{b}{c}
\right)^{(m-1)(p-1)/2}.
$$
Otherwise $(x_m/\phi_n)=1$ if $n$ has at least two distinct prime
factors. On the other hand if $n$ is even and $>2$ then  $(\phi_n/x_m)=1$, except when
$c\equiv -1 \pmod 4$, $n$ is a power of 2, and $m\equiv \pm 3
\pmod 8$, whence $(\phi_n/x_m)=-1$.
\endproclaim

Hence we can prove that $\phi_{p^k}$ is not a square  if $b$ is even and

$\bullet$\  $c\equiv 3 \pmod 4$, or

$\bullet$\  $4|b$ with $(b/c)=-1$ and $p\equiv 3 \pmod 4$, or

$\bullet$\  $b\equiv 2 \pmod 4$ with $(b/c)=-1$ and $p\equiv 7 \pmod 8$, or

$\bullet$\  $b\equiv 2 \pmod 4$  and $p\equiv 5 \pmod 8$, or

$\bullet$\  $b\equiv 2 \pmod 4$ with $(b/c)=1$ and $p\equiv 3 \pmod 8$.

\head 7. Open problems \endhead

We conjecture that for every non-periodic Lucas sequence $\{ x_n\}_{n\geq 0}$ there exists an integer $n_x$ such that if $n\geq n_{x}$ then $x_n$ has a primitive prime factor  that divides it to an odd power.  In Theorem 3 we proved this in the special case that $\Delta>0$ and $c\equiv 2 \pmod 4$, with $n_x=7$. Proposition 6.1 suggests that our approach is unlikely to yield the analogous result in all other cases where $2|bc$.  We were unable to give a  formula for the Jacobi symbol $(x_m/x_n)$ in general when $b$ and $c$ are odd (which includes the interesting case of the Fibonacci numbers) which can be used in this context (though see [6]), and we hope that others will embrace this challenge.

\enddocument